\documentclass[11pt]{amsart}
\usepackage{amsmath,amssymb,mathrsfs}
\newtheorem{theorem}{Theorem}
\newtheorem{proposition}[theorem]{Proposition}
\newtheorem{corollary}[theorem]{Corollary}

\newtheorem*{minoo}{Min-Oo's Conjecture}

\begin{document}

\title{Scalar curvature rigidity of geodesic balls in $S^n$}
\author{S. Brendle and F.C. Marques}
\address{Department of Mathematics \\ Stanford University \\ 450 Serra Mall, Bldg 380 \\ Stanford, CA 94305} 
\address{Instituto de Matem\'atica Pura e Aplicada (IMPA) \\ Estrada Dona Castorina 110 \\ 22460-320 Rio de Janeiro \\ Brazil}
\thanks{The first author was supported in part by the National Science Foundation under grant DMS-0905628. The second author was supported by CNPq-Brazil, FAPERJ, and the Stanford Department of Mathematics.}
\begin{abstract}
In this paper, we prove a scalar curvature rigidity result for geodesic balls in $S^n$. This result contrasts sharply with the recent counterexamples to Min-Oo's conjecture for the hemisphere (cf. \cite{Brendle-Marques-Neves}).
\end{abstract}
\maketitle

\section{Introduction}

In this paper, we study rigidity phenomena involving the scalar curvature. These questions are motivated to a large extent by the positive mass theorem in general relativity, which was proved by Schoen and Yau \cite{Schoen-Yau} and Witten \cite{Witten}. An important corollary of this theorem is that any Riemannian metric on $\mathbb{R}^n$ which has nonnegative scalar curvature and agrees with the Euclidean metric outside a compact set is necessarily flat.

It was observed by Miao \cite{Miao} that the positive mass theorem implies the following rigidity result for metrics on the unit ball:

\begin{theorem}
\label{rigidity.for.ball}
Suppose that $g$ is a Riemannian metric on the unit ball $B^n \subset \mathbb{R}^n$ with the following properties: 
\begin{itemize}
\item The scalar curvature of $g$ is nonnegative.
\item The induced metric on the boundary $\partial B^n$ agrees with the standard metric on $\partial B^n$.
\item The mean curvature of $\partial B^n$ with respect to $g$ is at least $n-1$.
\end{itemize}
Then $g$ is isometric to the standard metric on $B^n$.
\end{theorem}

An important generalization of Theorem \ref{rigidity.for.ball} was proved by Shi and Tam \cite{Shi-Tam}.

Motivated by the positive mass theorem, Min-Oo \cite{Min-Oo} posed the following question: 

\begin{minoo}
Suppose that $g$ is a Riemannian metric on the hemisphere $S_+^n$ with the following properties: 
\begin{itemize}
\item The scalar curvature of $g$ is at least $n(n-1)$.
\item The induced metric on the boundary $\partial S_+^n$ agrees with the standard metric on $\partial S_+^n$.
\item The boundary $\partial S_+^n$ is totally geodesic with respect to $g$.
\end{itemize}
Then $g$ is isometric to the standard metric on $S_+^n$.
\end{minoo}

Min-Oo's conjecture has been verified in many special cases (see e.g. \cite{Eichmair}, \cite{Hang-Wang1}, \cite{Hang-Wang2}). A related rigidity result for real projective space $\mathbb{RP}^3$ was established in \cite{Bray-Brendle-Eichmair-Neves}.

Very recently, counterexamples to Min-Oo's conjecture were constructed in \cite{Brendle-Marques-Neves}.

\begin{theorem}[S.~Brendle, F.C.~Marques, A.~Neves \cite{Brendle-Marques-Neves}]
\label{counterexample}
Given any integer $n \geq 3$, there exists a smooth Riemannian metric $\hat{g}$ on the hemisphere $S_+^n$ with the following properties:
\begin{itemize}
\item The scalar curvature of $\hat{g}$ is at least $n(n-1)$ at each point on $S_+^n$.
\item The scalar curvature of $\hat{g}$ is strictly greater than $n(n-1)$ at some point on $S_+^n$.
\item The metric $\hat{g}$ agrees with the standard metric in a neighborhood of $\partial S_+^n$.
\end{itemize}
\end{theorem}

The proof of Theorem \ref{counterexample} relies on a perturbation analysis.

In this paper, we study the analogous rigidity question for geodesic balls in $S^n$ of radius less than $\frac{\pi}{2}$. To fix notation, let $\overline{g}$ be the standard metric on $S^n$ and let $f: S^n \to \mathbb{R}$ denotes the restriction of the coordinate function $x_{n+1}$ to $S^n$. We will consider a domain of the form $\Omega = \{f \geq c\}$. If $c \geq \frac{2}{\sqrt{n+3}}$, we have the following rigidity result:

\begin{theorem}
\label{main.theorem}
Consider the domain $\Omega = \{f \geq c\}$, where $c \geq \frac{2}{\sqrt{n+3}}$. Let $g$ be a Riemannian metric on $\Omega$ with the following properties: 
\begin{itemize}
\item $R_g \geq n(n-1)$ at each point in $\Omega$. 
\item $H_g \geq H_{\overline{g}}$ at each point on $\partial \Omega$. 
\item The metrics $g$ and $\overline{g}$ induce the same metric on $\partial \Omega$. 
\end{itemize} 
If $g - \overline{g}$ is sufficiently small in the $C^2$-norm, then $\varphi^*(g) = \overline{g}$ for some diffeomorphism $\varphi: \Omega \to \Omega$ with $\varphi|_{\partial \Omega} = \text{\rm id}$.
\end{theorem}

We remark that the conclusion of Theorem \ref{main.theorem} holds under the weaker assumption that $g$ is close to $\overline{g}$ in $W^{2,p}$-norm for $p > n$. 

Note that Theorem \ref{main.theorem} is false for the hemisphere $\{f \geq 0\}$: by Theorem 4 in \cite{Brendle-Marques-Neves}, there exist Riemannian metrics on the hemisphere which satisfy the assumptions of Theorem \ref{main.theorem} and are arbitrary close to the standard metric $\overline{g}$ in the $C^\infty$-topology, but which are not isometric to $\overline{g}$.

The proof of Theorem \ref{main.theorem} relies on a perturbation analysis which is similar in spirit to Bartnik's work on the positive mass theorem (cf. \cite{Bartnik}, Section 5). Similar techniques have been employed in the study of the total scalar curvature functional (see e.g. \cite{Besse}, Section 4G) and the Yamabe flow (cf. \cite{Brendle}). Dai, Wang, and Wei \cite{Dai-Wang-Wei1},\cite{Dai-Wang-Wei2} have obtained interesting stability results for manifolds with parallel spinors, as well as for K\"ahler-Einstein manifolds.

\section{The scalar curvature and boundary mean curvature of a perturbed metric}

In this section, we consider a smooth manifold $\Omega$ with boundary $\partial \Omega$. Let $\overline{g}$ be a fixed Riemannian metric on $\Omega$. Moreover, we consider another Riemannian metric $g = \overline{g} + h$, where $|h|_{\overline{g}} \leq \frac{1}{2}$ at each point in $\Omega$. For abbreviation, we write $(h^2)_{ik} = \overline{g}^{jl} \, h_{ij} \, h_{kl}$.

\begin{proposition}
\label{scalar.curvature}
The scalar curvature of $g$ satisfies the pointwise estimate 
\begin{align*} 
&\Big | R_g - R_{\overline{g}} + \langle \text{\rm Ric}_{\overline{g}},h \rangle - \langle \text{\rm Ric}_{\overline{g}},h^2 \rangle \\ 
&+ \frac{1}{4} \, \overline{g}^{ij} \, \overline{g}^{kl} \, \overline{g}^{pq} \, \overline{D}_i h_{kp} \, \overline{D}_j h_{lq} - \frac{1}{2} \, \overline{g}^{ij} \, \overline{g}^{kl} \, \overline{g}^{pq} \, \overline{D}_i h_{kp} \, \overline{D}_l h_{jq} \\ 
&+ \frac{1}{4} \, \overline{g}^{pq} \, \partial_p(\text{\rm tr}_{\overline{g}}(h)) \, \partial_q(\text{\rm tr}_{\overline{g}}(h)) + \overline{D}_i \Big ( g^{ik} \, g^{jl} \, (\overline{D}_k h_{jl} - \overline{D}_l h_{jk}) \Big ) \Big | \\ 
&\leq C \, |h| \, |\overline{D} h|^2 + C \, |h|^3. 
\end{align*} 
Here, $\overline{D}$ denotes the Levi-Civita connection with respect to $\overline{g}$, and $C$ is a positive constant which depends only on $n$. 
\end{proposition}

\textbf{Proof.} 
The Levi-Civita connection with respect to $g$ is given by 
\[D_X Y = \overline{D}_X Y + \Gamma(X,Y),\] 
where $\Gamma$ is defined by 
\[g(\Gamma(X,Y),Z) = \frac{1}{2} \, \big ( (\overline{D}_X h)(Y,Z) + (\overline{D}_Y h)(X,Z) - (\overline{D}_Z h)(X,Y) \big ).\] 
In local coordinates, the tensor $\Gamma$ can be written in the form 
\[\Gamma_{jk}^m = \frac{1}{2} \, g^{lm} \, (\overline{D}_j h_{kl} + \overline{D}_k h_{jl} - \overline{D}_l h_{jk}).\] 
With this understood, the scalar curvature of $g$ is given by
\begin{align*} 
R_g &= g^{ik} \, (\text{\rm Ric}_{\overline{g}})_{ik} + g^{ik} \, g^{jl} \, g_{pq} \, \Gamma_{il}^q \, \Gamma_{jk}^p - g^{ik} \, g^{jl} \, g_{pq} \, \Gamma_{jl}^q \, \Gamma_{ik}^p \\ 
&- g^{ik} \, g^{jl} \, (\overline{D}_{i,k}^2 h_{jl} - \overline{D}_{i,l}^2 h_{jk}) 
\end{align*} 
(cf. \cite{Brendle-Marques-Neves}, Proposition 16). This implies 
\begin{align*} 
&\Big | R_g - R_{\overline{g}} + \langle \text{\rm Ric}_{\overline{g}},h \rangle - \langle \text{\rm Ric}_{\overline{g}},h^2 \rangle \\ 
&- \frac{3}{4} \, \overline{g}^{ij} \, \overline{g}^{kl} \, \overline{g}^{pq} \, \overline{D}_i h_{kp} \, \overline{D}_j h_{lq} + \frac{1}{2} \, \overline{g}^{ij} \, \overline{g}^{kl} \, \overline{g}^{pq} \, \overline{D}_i h_{kp} \, \overline{D}_l h_{jq} \\ 
&+ \frac{1}{4} \, \overline{g}^{pq} \, \partial_p(\text{\rm tr}_{\overline{g}}(h)) \, \partial_q(\text{\rm tr}_{\overline{g}}(h)) - \overline{g}^{ij} \, \overline{g}^{pq} \, \overline{D}_i h_{jp} \, \partial_q(\text{\rm tr}_{\overline{g}}(h)) \\ 
&+ \overline{g}^{ij} \, \overline{g}^{kl} \, \overline{g}^{pq} \, \overline{D}_i h_{jp} \, \overline{D}_k h_{lq} + g^{ik} \, g^{jl} \, (\overline{D}_{i,k}^2 h_{jl} - \overline{D}_{i,l}^2 h_{jk}) \Big | \\ 
&\leq C \, |h| \, |\overline{D} h|^2 + C \, |h|^3. 
\end{align*} 
Hence, we obtain 
\begin{align*} 
&\Big | R_g - R_{\overline{g}} + \langle \text{\rm Ric}_{\overline{g}},h \rangle - \langle \text{\rm Ric}_{\overline{g}},h^2 \rangle \\ 
&+ \frac{1}{4} \, \overline{g}^{ij} \, \overline{g}^{kl} \, \overline{g}^{pq} \, \overline{D}_i h_{kp} \, \overline{D}_j h_{lq} - \frac{1}{2} \, \overline{g}^{ij} \, \overline{g}^{kl} \, \overline{g}^{pq} \, \overline{D}_i h_{kp} \, \overline{D}_l h_{jq} \\ 
&+ \frac{1}{4} \, \overline{g}^{pq} \, \partial_p(\text{\rm tr}_{\overline{g}}(h)) \, \partial_q(\text{\rm tr}_{\overline{g}}(h)) + \overline{D}_i \Big ( g^{ik} \, g^{jl} \, (\overline{D}_k h_{jl} - \overline{D}_l h_{jk}) \Big ) \Big | \\ 
&\leq C \, |h| \, |\overline{D} h|^2 + C \, |h|^3, 
\end{align*} 
as claimed. \\

In the next step, we estimate the mean curvature of the boundary $\partial \Omega$ with respect to the metric $\overline{g}$. To that end, we assume that $g$ and $\overline{g}$ induce the same metric on the boundary $\partial \Omega$; in other words, we assume that $h(X,Y) = 0$ whenever $X$ and $Y$ are tangent vectors to $\partial \Omega$.

\begin{proposition}
\label{mean.curvature}
Assume that $g$ and $\overline{g}$ induce the same metric on the boundary $\partial \Omega$. Then the mean curvature of $\partial \Omega$ with respect to $g$ satisfies 
\begin{align*} 
&\bigg | 2 \, (H_g - H_{\overline{g}}) - \Big ( h(\overline{\nu},\overline{\nu}) - \frac{1}{4} \, h(\overline{\nu},\overline{\nu})^2 + \sum_{a=1}^{n-1} h(e_a,\overline{\nu})^2 \Big ) \, H_{\overline{g}} \\ 
&+ \big ( 1 - \frac{1}{2} \, h(\overline{\nu},\overline{\nu}) \big ) \, \sum_{a=1}^{n-1} \big ( 2 \, (\overline{D}_{e_a} h)(e_a,\overline{\nu}) - (\overline{D}_{\overline{\nu}} h)(e_a,e_a) \big ) \bigg | \\ 
&\leq C \, |h|^2 \, |\overline{D} h| + C \, |h|^3. 
\end{align*}
Here, $\{e_a: 1 \leq a \leq n-1\}$ is a local orthonormal frame on $\partial \Omega$, and $C$ is a positive constant that depends only on $n$.
\end{proposition}

\textbf{Proof.} 
Using the identity 
\[H_g \, \nu - H_{\overline{g}} \, \overline{\nu} = -\sum_{a=1}^{n-1} (D_{e_a} e_a - \overline{D}_{e_a} e_a) = -\sum_{a=1}^{n-1} \Gamma(e_a,e_a),\] 
we obtain 
\begin{align*} 
&2 \, \big ( H_g \, g(\nu,\overline{\nu}) - H_{\overline{g}} \, g(\overline{\nu},\overline{\nu}) \big ) \\ 
&= -2 \, \sum_{a=1}^{n-1} g(\Gamma(e_a,e_a),\overline{\nu}) = -\sum_{a=1}^{n-1} \big ( 2 \,(\overline{D}_{e_a} h)(e_a,\overline{\nu}) - (\overline{D}_{\overline{\nu}} h)(e_a,e_a) \big ). 
\end{align*}
Clearly, $g(\overline{\nu},\overline{\nu}) = 1 + h(\overline{\nu},\overline{\nu})$. Moreover, it is easy to see that the vector $\overline{\nu} - \sum_{a=1}^{n-1} h(e_a,\overline{\nu}) \, e_a$ is orthogonal to $\partial \Omega$ with respect to $g$. From this, we deduce that 
\[\overline{\nu} - \sum_{a=1}^{n-1} h(e_a,\overline{\nu}) \, e_a = \bigg ( 1 + h(\overline{\nu},\overline{\nu}) - \sum_{a=1}^{n-1} h(e_a,\overline{\nu})^2 \bigg )^{\frac{1}{2}} \, \nu,\] 
hence 
\[g(\nu,\overline{\nu}) = \bigg ( 1 + h(\overline{\nu},\overline{\nu}) - \sum_{a=1}^{n-1} h(e_a,\overline{\nu})^2 \bigg )^{\frac{1}{2}}.\] 
Substituting these identities into the previous formula for $H_g$, the assertion follows.

\section{Perturbations of the standard metric on $S^n$}

We now consider perturbations of the standard metric $\overline{g}$ on $S^n$. To fix notation, let $f: S^n \to \mathbb{R}$ denote the restriction of the coordinate function $x_{n+1}$ to $S^n$, and let $\Omega = \{f \geq c\}$ be a geodesic ball centered at the north pole. Here, $c$ is a positive real number which will be specified later. 

Let $g$ be a Riemannian metric on $\Omega$. We will assume throughout that $g$ and $\overline{g}$ induce the same metric on the boundary $\partial \Omega$. Moreover, we assume that $g = \overline{g} + h$, where $|h|_{\overline{g}} \leq \frac{1}{2}$ at each point in $\Omega$. 

Our goal in this section is to estimate the integral 
\[\int_\Omega (R_g - n(n-1)) \, f \, d\text{\rm vol}_{\overline{g}}\] 
(see also \cite{Fischer-Marsden}).

\begin{proposition}
\label{scalar.curvature.2}
We have 
\begin{align*} 
&\bigg | \int_\Omega (R_g - n(n-1) - (n-1) \, |h|_{\overline{g}}^2) \, f \, d\text{\rm vol}_{\overline{g}} + \frac{1}{4} \int_\Omega |\overline{D} h|^2 \, f \, d\text{\rm vol}_{\overline{g}} \\ 
&+ \frac{1}{4} \int_\Omega |\overline{\nabla} (\text{\rm tr}_{\overline{g}}(h))|^2 \, f \, d\text{\rm vol}_{\overline{g}} - \frac{1}{2} \int_\Omega \overline{g}^{ij} \, \overline{g}^{kl} \, \overline{g}^{pq} \, \overline{D}_i h_{kp} \, \overline{D}_l h_{jq}  \, f \, d\text{\rm vol}_{\overline{g}} \\ 
&+ \int_\Omega \overline{g}^{ik} \, \overline{g}^{jp} \, \overline{g}^{lq} \, h_{pq} \, (\overline{D}_k h_{jl} - \overline{D}_l h_{jk}) \, \partial_i f \, d\text{\rm vol}_{\overline{g}} \\ 
&+ \int_\Omega \overline{g}^{ip} \, \overline{g}^{kq} \, \overline{g}^{jl} \, h_{pq} \, (\overline{D}_k h_{jl} - \overline{D}_l h_{jk}) \, \partial_i f \, d\text{\rm vol}_{\overline{g}} \\ 
&+ \int_{\partial \Omega} \overline{g}^{jl} \, (\overline{D}_k h_{jl} - \overline{D}_l h_{jk}) \, \overline{\nu}^k \, f \, d\sigma_{\overline{g}} \\ 
&- \int_{\partial \Omega} \overline{g}^{jp} \, \overline{g}^{lq} \, h_{pq} \, (\overline{D}_k h_{jl} - \overline{D}_l h_{jk}) \, \overline{\nu}^k \, f \, d\sigma_{\overline{g}} \\ 
&- \int_{\partial \Omega} \overline{g}^{kq} \, \overline{g}^{jl} \, h_{pq} \, (\overline{D}_k h_{jl} - \overline{D}_l h_{jk}) \, \overline{\nu}^p \, f \, d\sigma_{\overline{g}} \bigg | \\ 
&\leq C \int_\Omega |h| \, |\overline{D} h|^2 \, d\text{\rm vol}_{\overline{g}} + C \int_\Omega |h|^3 \, d\text{\rm vol}_{\overline{g}} + C \int_{\partial \Omega} |h|^2 \, |\overline{D} h| \, d\sigma_{\overline{g}}, 
\end{align*} 
where $C$ is a positive constant that depends only on $n$ and $c$.
\end{proposition}

\textbf{Proof.} Using Proposition \ref{scalar.curvature} and the divergence theorem, we obtain 
\begin{align*} 
&\bigg | \int_\Omega (R_g - n(n-1) + (n-1) \, \text{\rm tr}_{\overline{g}}(h) - (n-1) \, |h|_{\overline{g}}^2) \, f \, d\text{\rm vol}_{\overline{g}} \\ 
&+ \frac{1}{4} \int_\Omega |\overline{D} h|^2 \, f \, d\text{\rm vol}_{\overline{g}} - \frac{1}{2} \int_\Omega \overline{g}^{ij} \, \overline{g}^{kl} \, \overline{g}^{pq} \, \overline{D}_i h_{kp} \, \overline{D}_l h_{jq}  \, f \, d\text{\rm vol}_{\overline{g}} \\ 
&+ \frac{1}{4} \int_\Omega |\overline{\nabla} (\text{\rm tr}_{\overline{g}}(h))|^2 \, f \, d\text{\rm vol}_{\overline{g}} - \int_\Omega g^{ik} \, g^{jl} \, (\overline{D}_k h_{jl} - \overline{D}_l h_{jk}) \, \partial_i f \, d\text{\rm vol}_{\overline{g}} \\ 
&+ \int_{\partial \Omega} g^{ik} \, g^{jl} \, (\overline{D}_k h_{jl} - \overline{D}_l h_{jk}) \, \overline{g}_{im} \, \overline{\nu}^m \, f \, d\sigma_{\overline{g}} \bigg | \\ 
&\leq C \int_\Omega |h| \, |\overline{D} h|^2 \, d\text{\rm vol}_{\overline{g}} + C \int_\Omega |h|^3 \, d\text{\rm vol}_{\overline{g}}. 
\end{align*} 
Here, $\overline{\nu}$ denotes the outward-pointing unit normal vector to $\partial \Omega$ with respect to the metric $\overline{g}$. Using the identity $\overline{D}_{i,k}^2 f = -f \, \overline{g}_{ik}$, we obtain 
\begin{align*} 
&\int_\Omega (n-1) \, \text{\rm tr}_{\overline{g}}(h) \, f \, d\text{\rm vol}_{\overline{g}} - \int_\Omega \overline{g}^{ik} \, \overline{g}^{jl} \, (\overline{D}_k h_{jl} - \overline{D}_l h_{jk}) \, \partial_i f \, d\text{\rm vol}_{\overline{g}} \\ 
&= -\int_{\partial \Omega} (\text{\rm tr}_{\overline{g}}(h) \, \partial_{\overline{\nu}} f - h(\overline{\nu},\overline{\nabla} f)) \, d\sigma_{\overline{g}} = 0. 
\end{align*}
Thus, we conclude that 
\begin{align*} 
&\bigg | \int_\Omega (R_g - n(n-1) - (n-1) \, |h|_{\overline{g}}^2) \, f \, d\text{\rm vol}_{\overline{g}} + \frac{1}{4} \int_\Omega |\overline{D} h|^2 \, f \, d\text{\rm vol}_{\overline{g}} \\ 
&+ \frac{1}{4} \int_\Omega |\overline{\nabla} (\text{\rm tr}_{\overline{g}}(h))|^2 \, f \, d\text{\rm vol}_{\overline{g}} - \frac{1}{2} \int_\Omega \overline{g}^{ij} \, \overline{g}^{kl} \, \overline{g}^{pq} \, \overline{D}_i h_{kp} \, \overline{D}_l h_{jq}  \, f \, d\text{\rm vol}_{\overline{g}} \\ 
&- \int_\Omega (g^{ik} \, g^{jl} - \overline{g}^{ik} \, \overline{g}^{jl}) \, (\overline{D}_k h_{jl} - \overline{D}_l h_{jk}) \, \partial_i f \, d\text{\rm vol}_{\overline{g}} \\ 
&+ \int_{\partial \Omega} g^{ik} \, g^{jl} \, (\overline{D}_k h_{jl} - \overline{D}_l h_{jk}) \, \overline{g}_{im} \, \overline{\nu}^m \, f \, d\sigma_{\overline{g}} \bigg | \\ 
&\leq C \int_\Omega |h| \, |\overline{D} h|^2 \, d\text{\rm vol}_{\overline{g}} + C \int_\Omega |h|^3 \, d\text{\rm vol}_{\overline{g}}. 
\end{align*} From this, the assertion follows easily. \\

In the remainder of this section, we will assume that $h$ is divergence-free in the sense that $\overline{g}^{ik} \, \overline{D}_i h_{kl} = 0$. 

\begin{proposition}
\label{scalar.curvature.3}
Assume that $h$ is divergence-free. Then 
\begin{align*} 
&\bigg | \int_\Omega (R_g - n(n-1)) \, f \, d\text{\rm vol}_{\overline{g}} + \frac{1}{4} \int_\Omega |\overline{D} h|^2 \, f \, d\text{\rm vol}_{\overline{g}} \\ 
&+ \frac{1}{4} \int_\Omega |\overline{\nabla} (\text{\rm tr}_{\overline{g}}(h))|^2 \, f \, d\text{\rm vol}_{\overline{g}} + \frac{1}{2} \int_\Omega |h|_{\overline{g}}^2 \, f \, d\text{\rm vol}_{\overline{g}} \\ 
&+ \frac{1}{2} \int_\Omega \text{\rm tr}_{\overline{g}}(h)^2 \, f \, d\text{\rm vol}_{\overline{g}} + \frac{1}{4} \int_{\partial \Omega} (|h|_{\overline{g}}^2 + 3 \, h(\overline{\nu},\overline{\nu})^2) \, \partial_{\overline{\nu}} f \, d\sigma_{\overline{g}} \\ 
&+ \int_{\partial \Omega} \overline{g}^{jl} \, \overline{D}_k h_{jl} \, \overline{\nu}^k \, f \, d\sigma_{\overline{g}} - \frac{1}{2} \int_{\partial \Omega} \overline{g}^{kl} \, \overline{g}^{pq} \, h_{kp} \, \overline{D}_l h_{jq} \, \overline{\nu}^j \, f \, d\sigma_{\overline{g}} \\ 
&- \int_{\partial \Omega} \overline{g}^{jp} \, \overline{g}^{lq} \, h_{pq} \, (\overline{D}_k h_{jl} - \overline{D}_l h_{jk}) \, \overline{\nu}^k \, f \, d\sigma_{\overline{g}} \\ 
&- \int_{\partial \Omega} \overline{g}^{kq} \, \overline{g}^{jl} \, h_{pq} \, \overline{D}_k h_{jl} \, \overline{\nu}^p \, f \, d\sigma_{\overline{g}} \bigg | \\ 
&\leq C \int_\Omega |h| \, |\overline{D} h|^2 \, d\text{\rm vol}_{\overline{g}} + C \int_\Omega |h|^3 \, d\text{\rm vol}_{\overline{g}} + C \int_{\partial \Omega} |h|^2 \, |\overline{D} h| \, d\sigma_{\overline{g}}, 
\end{align*} 
where $C$ is a positive constant that depends only on $n$ and $c$.
\end{proposition}

\textbf{Proof.} 
Since $\overline{g}$ has constant sectional curvature $1$, we have 
\[\overline{D}_{i,l}^2 h_{jq} = \overline{D}_{l,i}^2 h_{jq} + h_{lq} \, \overline{g}_{ij} - h_{iq} \, \overline{g}_{jl} + h_{jl} \, \overline{g}_{iq} - h_{ij} \, \overline{g}_{lq}.\]
Since $h$ is divergence-free, it follows that  
\[\overline{g}^{ij} \, \overline{D}_{i,l}^2 h_{jq} = n \, h_{lq} - \text{\rm tr}_{\overline{g}}(h) \, \overline{g}_{lq}.\] 
This implies 
\begin{align*} 
&-\int_\Omega \overline{g}^{ij} \, \overline{g}^{kl} \, \overline{g}^{pq} \, h_{kp} \, \overline{D}_l h_{jq} \, \partial_i f \, d\text{\rm vol}_{\overline{g}} - \int_\Omega \overline{g}^{ij} \, \overline{g}^{kl} \, \overline{g}^{pq} \, \overline{D}_i h_{kp} \, \overline{D}_l h_{jq} \, f \, d\text{\rm vol}_{\overline{g}} \\ 
&= \int_\Omega \overline{g}^{ij} \, \overline{g}^{kl} \, \overline{g}^{pq} \, h_{kp} \, \overline{D}_{i,l}^2 h_{jq} \, f \, d\text{\rm vol}_{\overline{g}} - \int_{\partial \Omega} \overline{g}^{kl} \, \overline{g}^{pq} \, h_{kp} \, \overline{D}_l h_{jq} \, \overline{\nu}^j \, f \, d\sigma_{\overline{g}} \\ 
&= n \int_\Omega |h|_{\overline{g}}^2 \, f \, d\text{\rm vol}_{\overline{g}} - \int_\Omega \text{\rm tr}_{\overline{g}}(h)^2 \, f \, d\text{\rm vol}_{\overline{g}} - \int_{\partial \Omega} \overline{g}^{kl} \, \overline{g}^{pq} \, h_{kp} \, \overline{D}_l h_{jq} \, \overline{\nu}^j \, f \, d\sigma_{\overline{g}}. 
\end{align*} From this, we deduce that 
\begin{align*} 
&\int_\Omega \overline{g}^{ik} \, \overline{g}^{jp} \, \overline{g}^{lq} \, h_{pq} \, (\overline{D}_k h_{jl} - \overline{D}_l h_{jk}) \, \partial_i f \, d\text{\rm vol}_{\overline{g}} \\ 
&- \frac{1}{2} \int_\Omega \overline{g}^{ij} \, \overline{g}^{kl} \, \overline{g}^{pq} \, \overline{D}_i h_{kp} \, \overline{D}_l h_{jq} \, f \, d\text{\rm vol}_{\overline{g}} \\ 
&= \frac{1}{2} \int_\Omega \overline{g}^{ik} \, \partial_k(|h|_{\overline{g}}^2) \, \partial_i f \, d\text{\rm vol}_{\overline{g}} - \frac{1}{2} \int_\Omega \overline{g}^{ik} \, \overline{g}^{jp} \, \overline{g}^{lq} \, h_{pq} \, \overline{D}_l h_{jk} \, \partial_i f \, d\text{\rm vol}_{\overline{g}} \\ 
&+ \frac{n}{2} \int_\Omega |h|_{\overline{g}}^2 \, f \, d\text{\rm vol}_{\overline{g}} - \frac{1}{2} \int_\Omega \text{\rm tr}_{\overline{g}}(h)^2 \, f \, d\text{\rm vol}_{\overline{g}} \\ 
&- \frac{1}{2} \int_{\partial \Omega} \overline{g}^{kl} \, \overline{g}^{pq} \, h_{kp} \, \overline{D}_l h_{jq} \, \overline{\nu}^j \, f \, d\sigma_{\overline{g}}. 
\end{align*} 
Integration by parts gives 
\begin{align*} 
&\int_\Omega \overline{g}^{ik} \, \overline{g}^{jp} \, \overline{g}^{lq} \, h_{pq} \, (\overline{D}_k h_{jl} - \overline{D}_l h_{jk}) \, \partial_i f \, d\text{\rm vol}_{\overline{g}} \\ 
&- \frac{1}{2} \int_\Omega \overline{g}^{ij} \, \overline{g}^{kl} \, \overline{g}^{pq} \, \overline{D}_i h_{kp} \, \overline{D}_l h_{jq} \, f \, d\text{\rm vol}_{\overline{g}} \\ 
&= -\frac{1}{2} \int_\Omega |h|_{\overline{g}}^2 \, \Delta_{\overline{g}} f \, d\text{\rm vol}_{\overline{g}} + \frac{1}{2} \int_\Omega \overline{g}^{ik} \, \overline{g}^{jp} \, \overline{g}^{lq} \, h_{pq} \, h_{jk} \, \overline{D}_{i,l}^2 f \, d\text{\rm vol}_{\overline{g}} \\ 
&+ \frac{1}{2} \int_{\partial \Omega} |h|_{\overline{g}}^2 \, \partial_{\overline{\nu}} f \, d\sigma_{\overline{g}} - \frac{1}{2} \int_{\partial \Omega} \overline{g}^{ik} \, \overline{g}^{jp} \, h_{pq} \, h_{jk} \, \partial_i f \, \overline{\nu}^q \, d\sigma_{\overline{g}} \\ 
&+ \frac{n}{2} \int_\Omega |h|_{\overline{g}}^2 \, f \, d\text{\rm vol}_{\overline{g}} - \frac{1}{2} \int_\Omega \text{\rm tr}_{\overline{g}}(h)^2 \, f \, d\text{\rm vol}_{\overline{g}} \\ 
&- \frac{1}{2} \int_{\partial \Omega} \overline{g}^{kl} \, \overline{g}^{pq} \, h_{kp} \, \overline{D}_l h_{jq} \, \overline{\nu}^j \, f \, d\sigma_{\overline{g}} \\ 
&= \frac{2n-1}{2} \int_\Omega |h|_{\overline{g}}^2 \, f \, d\text{\rm vol}_{\overline{g}} - \frac{1}{2} \int_\Omega \text{\rm tr}_{\overline{g}}(h)^2 \, f \, d\text{\rm vol}_{\overline{g}} \\ 
&+ \frac{1}{2} \int_{\partial \Omega} |h|_{\overline{g}}^2 \, \partial_{\overline{\nu}} f \, d\sigma_{\overline{g}} - \frac{1}{2} \int_{\partial \Omega} \overline{g}^{ik} \, \overline{g}^{jp} \, h_{pq} \, h_{jk} \, \partial_i f \, \overline{\nu}^q \, d\sigma_{\overline{g}} \\ 
&- \frac{1}{2} \int_{\partial \Omega} \overline{g}^{kl} \, \overline{g}^{pq} \, h_{kp} \, \overline{D}_l h_{jq} \, \overline{\nu}^j \, f \, d\sigma_{\overline{g}}. 
\end{align*}
Moreover, we have 
\begin{align*} 
&\int_\Omega \overline{g}^{ip} \, \overline{g}^{kq} \, \overline{g}^{jl} \, h_{pq} \, (\overline{D}_k h_{jl} - \overline{D}_l h_{jk}) \, \partial_i f \, d\text{\rm vol}_{\overline{g}} \\ 
&= \int_\Omega \overline{g}^{ip} \, \overline{g}^{kq} \, h_{pq} \, \partial_k(\text{\rm tr}_{\overline{g}}(h)) \, \partial_i f \, d\text{\rm vol}_{\overline{g}} \\ 
&= -\int_\Omega \overline{g}^{ip} \, \overline{g}^{kq} \, h_{pq} \, \text{\rm tr}_{\overline{g}}(h) \, \overline{D}_{i,k}^2 f \, d\text{\rm vol}_{\overline{g}} + \int_{\partial \Omega} \overline{g}^{ip} \, h_{pq} \, \text{\rm tr}_{\overline{g}}(h) \, \partial_i f \, \overline{\nu}^q \, d\sigma_{\overline{g}} \\ 
&= \int_\Omega \text{\rm tr}_{\overline{g}}(h)^2 \, f \, d\text{\rm vol}_{\overline{g}} + \int_{\partial \Omega} \overline{g}^{ip} \, h_{pq} \, \text{\rm tr}_{\overline{g}}(h) \, \partial_i f \, \overline{\nu}^q \, d\sigma_{\overline{g}}. 
\end{align*} 
Putting these facts together, we obtain 
\begin{align*} 
&\int_\Omega \overline{g}^{ip} \, \overline{g}^{kq} \, \overline{g}^{jl} \, h_{pq} \, (\overline{D}_k h_{jl} - \overline{D}_l h_{jk}) \, \partial_i f \, d\text{\rm vol}_{\overline{g}} \\ 
&+ \int_\Omega \overline{g}^{ik} \, \overline{g}^{jp} \, \overline{g}^{lq} \, h_{pq} \, (\overline{D}_k h_{jl} - \overline{D}_l h_{jk}) \, \partial_i f \, d\text{\rm vol}_{\overline{g}} \\ 
&- \frac{1}{2} \int_\Omega \overline{g}^{ij} \, \overline{g}^{kl} \, \overline{g}^{pq} \, \overline{D}_i h_{kp} \, \overline{D}_l h_{jq} \, f \, d\text{\rm vol}_{\overline{g}} \\ 
&= \frac{2n-1}{2} \int_\Omega |h|_{\overline{g}}^2 \, f \, d\text{\rm vol}_{\overline{g}} + \frac{1}{2} \int_\Omega \text{\rm tr}_{\overline{g}}(h)^2 \, f \, d\text{\rm vol}_{\overline{g}} \\ 
&+ \frac{1}{2} \int_{\partial \Omega} |h|_{\overline{g}}^2 \, \partial_{\overline{\nu}} f \, d\sigma_{\overline{g}} - \frac{1}{2} \int_{\partial \Omega} \overline{g}^{ik} \, \overline{g}^{jp} \, h_{pq} \, h_{jk} \, \partial_i f \, \overline{\nu}^q \, d\sigma_{\overline{g}} \\ 
&+ \int_{\partial \Omega} \overline{g}^{ip} \, h_{pq} \, \text{\rm tr}_{\overline{g}}(h) \, \partial_i f \, \overline{\nu}^q \, d\sigma_{\overline{g}} - \frac{1}{2} \int_{\partial \Omega} \overline{g}^{kl} \, \overline{g}^{pq} \, h_{kp} \, \overline{D}_l h_{jq} \, \overline{\nu}^j \, f \, d\sigma_{\overline{g}} \\ 
&= \frac{2n-1}{2} \int_\Omega |h|_{\overline{g}}^2 \, f \, d\text{\rm vol}_{\overline{g}} + \frac{1}{2} \int_\Omega \text{\rm tr}_{\overline{g}}(h)^2 \, f \, d\text{\rm vol}_{\overline{g}} \\ 
&+ \frac{1}{4} \int_{\partial \Omega} (|h|_{\overline{g}}^2 + 3 \, h(\overline{\nu},\overline{\nu})^2) \, \partial_{\overline{\nu}} f \, d\sigma_{\overline{g}} - \frac{1}{2} \int_{\partial \Omega} \overline{g}^{kl} \, \overline{g}^{pq} \, h_{kp} \, \overline{D}_l h_{jq} \, \overline{\nu}^j \, f \, d\sigma_{\overline{g}}. 
\end{align*} 
Hence, the assertion follows from Proposition \ref{scalar.curvature.2}.

\section{Analysis of the boundary terms}

In this section, we analyze the boundary terms in Proposition \ref{scalar.curvature.3}. As in the previous section, we assume that $\overline{g}$ is the standard metric on $S^n$, and $\Omega = \{f \geq c\}$ centered at the north pole. Moreover, we consider a Riemannian metric on $\Omega$ of the form $g = \overline{g} + h$, where $|h|_{\overline{g}} \leq \frac{1}{2}$ at each point in $\Omega$. 

\begin{proposition}
\label{boundary.terms.1}
Assume that $h$ is divergence-free. Then 
\begin{align*} 
&\int_{\partial \Omega} \overline{g}^{jl} \, \overline{D}_k h_{jl} \, \overline{\nu}^k \, d\sigma_{\overline{g}} - \frac{1}{2} \int_{\partial \Omega} \overline{g}^{kl} \, \overline{g}^{pq} \, h_{kp} \, \overline{D}_l h_{jq} \, \overline{\nu}^j \, d\sigma_{\overline{g}} \\ 
&- \int_{\partial \Omega} \overline{g}^{jp} \, \overline{g}^{lq} \, h_{pq} \, (\overline{D}_k h_{jl} - \overline{D}_l h_{jk}) \, \overline{\nu}^k \, d\sigma_{\overline{g}} \\ 
&- \int_{\partial \Omega} \overline{g}^{kq} \, \overline{g}^{jl} \, h_{pq} \, \overline{D}_k h_{jl} \, \overline{\nu}^p \, d\sigma_{\overline{g}} \\ 
&= - \int_{\partial \Omega} (1 - h(\overline{\nu},\overline{\nu})) \, \sum_{a=1}^{n-1} \big ( 2 \, (\overline{D}_{e_a} h)(e_a,\overline{\nu}) - (\overline{D}_{\overline{\nu}} h)(e_a,e_a) \big ) \, d\sigma_{\overline{g}} \\ 
&+ \int_{\partial \Omega} \big ( 1 - \frac{1}{2} \, h(\overline{\nu},\overline{\nu}) \big ) \, h(\overline{\nu},\overline{\nu}) \, H_{\overline{g}} \, d\sigma_{\overline{g}} \\ 
&+ \frac{3n-2}{2(n-1)} \int_{\partial \Omega} \sum_{a=1}^{n-1} h(e_a,\overline{\nu})^2 \, H_{\overline{g}} \, d\sigma_{\overline{g}}. 
\end{align*} 
Here, $\{e_a: 1 \leq a \leq n-1\}$ is a local orthonormal frame on $\partial \Omega$, and $C$ is a positive constant that depends only on $n$ and $c$.
\end{proposition}

\textbf{Proof.} 
Let $\{e_a: 1 \leq a \leq n-1\}$ be a local orthonormal frame on $\partial \Omega$. Since $h$ is divergence-free, we have 
\begin{align*}
&\overline{g}^{jl} \, \overline{D}_k h_{jl} \, \overline{\nu}^k - \frac{1}{2} \, \overline{g}^{kl} \, \overline{g}^{pq} \, h_{kp} \, \overline{D}_l h_{jq} \, \overline{\nu}^j \\ 
&- \overline{g}^{jp} \, \overline{g}^{lq} \, h_{pq} \, (\overline{D}_k h_{jl} - \overline{D}_l h_{jk}) \, \overline{\nu}^k - \overline{g}^{kq} \, \overline{g}^{jl} \, h_{pq} \, \overline{D}_k h_{jl} \, \overline{\nu}^p \\ 
&= -(1 - h(\overline{\nu},\overline{\nu})) \, \sum_{a=1}^{n-1} \big ( 2 \, (\overline{D}_{e_a} h)(e_a,\overline{\nu}) - (\overline{D}_{\overline{\nu}} h)(e_a,e_a) \big ) \\ 
&+ \big ( 1 - \frac{1}{2} \, h(\overline{\nu},\overline{\nu}) \big ) \sum_{a=1}^{n-1} (\overline{D}_{e_a} h)(e_a,\overline{\nu}) - \frac{1}{2} \sum_{a=1}^{n-1} (\overline{D}_{e_a} h)(\overline{\nu},\overline{\nu}) \, h(e_a,\overline{\nu}) \\ 
&+ \frac{3}{2} \sum_{a,b=1}^{n-1} h(e_a,\overline{\nu}) \, (\overline{D}_{e_b} h)(e_a,e_b) - \sum_{a,b=1}^{n-1} h(e_a,\overline{\nu}) \, (\overline{D}_{e_a} h)(e_b,e_b). 
\end{align*} 
At this point, we define a one-form $\omega$ on $\partial \Omega$ by $\omega(e_a) = (1 - \frac{1}{2} \, h(\overline{\nu},\overline{\nu})) \, h(e_a,\overline{\nu})$. Since $\partial \Omega$ is umbilic with respect to $\overline{g}$, we have 
\[\overline{D}_{e_a} \overline{\nu} = \frac{1}{n-1} \, H_{\overline{g}} \, e_a,\] 
where $H_{\overline{g}}$ denotes the mean curvature of $\partial \Omega$ with respect to the metric $\overline{g}$. Using this relation, we obtain the following formula for the divergence of $\omega$: 
\begin{align*} 
\text{\rm div}_{\partial \Omega}(\omega) 
&= \big ( 1 - \frac{1}{2} \, h(\overline{\nu},\overline{\nu}) \big ) \sum_{a=1}^{n-1} (\overline{D}_{e_a} h)(e_a,\overline{\nu}) - \frac{1}{2} \sum_{a=1}^{n-1} (\overline{D}_{e_a} h)(\overline{\nu},\overline{\nu}) \, h(e_a,\overline{\nu}) \\ 
&- \big ( 1 - \frac{1}{2} \, h(\overline{\nu},\overline{\nu}) \big ) \, h(\overline{\nu},\overline{\nu}) \, H_{\overline{g}} - \frac{1}{n-1} \sum_{a=1}^{n-1} h(e_a,\overline{\nu})^2 \, H_{\overline{g}}. 
\end{align*} 
Moreover, we have the pointwise identities 
\[\sum_{b=1}^{n-1} (\overline{D}_{e_b} h)(e_a,e_b) = \frac{n}{n-1} \, h(e_a,\overline{\nu}) \, H_{\overline{g}}\] 
and 
\[\sum_{b=1}^{n-1} (\overline{D}_{e_a} h)(e_b,e_b) = \frac{2}{n-1} \, h(e_a,\overline{\nu}) \, H_{\overline{g}}.\]
Putting these facts together, we obtain 
\begin{align*}
&\overline{g}^{jl} \, \overline{D}_k h_{jl} \, \overline{\nu}^k - \frac{1}{2} \, \overline{g}^{kl} \, \overline{g}^{pq} \, h_{kp} \, \overline{D}_l h_{jq} \, \overline{\nu}^j \\ 
&- \overline{g}^{jp} \, \overline{g}^{lq} \, h_{pq} \, (\overline{D}_k h_{jl} - \overline{D}_l h_{jk}) \, \overline{\nu}^k - \overline{g}^{kq} \, \overline{g}^{jl} \, h_{pq} \, \overline{D}_k h_{jl} \, \overline{\nu}^p \\ 
&= -(1 - h(\overline{\nu},\overline{\nu})) \, \sum_{a=1}^{n-1} \big ( 2 \, (\overline{D}_{e_a} h)(e_a,\overline{\nu}) - (\overline{D}_{\overline{\nu}} h)(e_a,e_a) \big ) \\ 
&+ \big ( 1 - \frac{1}{2} \, h(\overline{\nu},\overline{\nu}) \big ) \, h(\overline{\nu},\overline{\nu}) \, H_{\overline{g}} + \frac{3n-2}{2(n-1)} \sum_{a=1}^{n-1} h(e_a,\overline{\nu})^2 \, H_{\overline{g}} + \text{\rm div}_{\partial \Omega}(\omega). 
\end{align*} 
Therefore, the assertion follows from the divergence theorem. \\

Combining Proposition \ref{boundary.terms.1} and Proposition \ref{mean.curvature}, we can draw the following conclusion:

\begin{corollary} 
\label{boundary.terms.2}
If $h$ is divergence-free, then we have 
\begin{align*} 
&\bigg | \int_{\partial \Omega} (2 - h(\overline{\nu},\overline{\nu})) \, (H_g - H_{\overline{g}}) \, d\sigma_{\overline{g}} \\ 
&- \int_{\partial \Omega} \overline{g}^{jl} \, \overline{D}_k h_{jl} \, \overline{\nu}^k \, d\sigma_{\overline{g}} + \frac{1}{2} \int_{\partial \Omega} \overline{g}^{kl} \, \overline{g}^{pq} \, h_{kp} \, \overline{D}_l h_{jq} \, \overline{\nu}^j \, d\sigma_{\overline{g}} \\ 
&+ \int_{\partial \Omega} \overline{g}^{jp} \, \overline{g}^{lq} \, h_{pq} \, (\overline{D}_k h_{jl} - \overline{D}_l h_{jk}) \, \overline{\nu}^k \, d\sigma_{\overline{g}} \\ 
&+ \int_{\partial \Omega} \overline{g}^{kq} \, \overline{g}^{jl} \, h_{pq} \, \overline{D}_k h_{jl} \, \overline{\nu}^p \, d\sigma_{\overline{g}} \\ 
&+ \frac{1}{4} \int_{\partial \Omega} h(\overline{\nu},\overline{\nu})^2 \, H_{\overline{g}} \, d\sigma_{\overline{g}} + \frac{n}{2(n-1)} \int_{\partial \Omega} \sum_{a=1}^{n-1} h(e_a,\overline{\nu})^2 \, H_{\overline{g}} \, d\sigma_{\overline{g}} \bigg | \\ 
&\leq C \int_{\partial \Omega} |h|^2 \, |\overline{D} h| \, d\sigma_{\overline{g}} + C \int_{\partial \Omega} |h|^3 \, d\sigma_{\overline{g}}, 
\end{align*} 
where $C$ is a positive constant that depends only on $n$ and $c$.
\end{corollary}

\textbf{Proof.} 
It follows from Proposition \ref{mean.curvature} that 
\begin{align*} 
&\bigg | \int_{\partial \Omega} (2 - h(\overline{\nu},\overline{\nu})) \, (H_g - H_{\overline{g}}) \, d\sigma_{\overline{g}} \\ 
&- \int_{\partial \Omega} \Big ( h(\overline{\nu},\overline{\nu}) - \frac{3}{4} \, h(\overline{\nu},\overline{\nu})^2 + \sum_{a=1}^{n-1} h(e_a,\overline{\nu})^2 \Big ) \, H_{\overline{g}} \, d\sigma_{\overline{g}} \\ 
&+ \int_{\partial \Omega} (1 - h(\overline{\nu},\overline{\nu})) \sum_{a=1}^{n-1} \big ( 2 \, (\overline{D}_{e_a} h)(e_a,\overline{\nu}) - (\overline{D}_{\overline{\nu}} h)(e_a,e_a) \big ) \, d\sigma_{\overline{g}} \bigg | \\ 
&\leq C \int_{\partial \Omega} |h|^2 \, |\overline{D} h| \, d\sigma_{\overline{g}} + C \int_{\partial \Omega} |h|^3 \, d\sigma_{\overline{g}}. 
\end{align*}
Moreover, we have 
\begin{align*} 
&\int_{\partial \Omega} \overline{g}^{jl} \, \overline{D}_k h_{jl} \, \overline{\nu}^k \, d\sigma_{\overline{g}} - \frac{1}{2} \int_{\partial \Omega} \overline{g}^{kl} \, \overline{g}^{pq} \, h_{kp} \, \overline{D}_l h_{jq} \, \overline{\nu}^j \, d\sigma_{\overline{g}} \\ 
&- \int_{\partial \Omega} \overline{g}^{jp} \, \overline{g}^{lq} \, h_{pq} \, (\overline{D}_k h_{jl} - \overline{D}_l h_{jk}) \, \overline{\nu}^k \, d\sigma_{\overline{g}} \\ 
&- \int_{\partial \Omega} \overline{g}^{kq} \, \overline{g}^{jl} \, h_{pq} \, \overline{D}_k h_{jl} \, \overline{\nu}^p \, d\sigma_{\overline{g}} \\ 
&= - \int_{\partial \Omega} (1 - h(\overline{\nu},\overline{\nu})) \, \sum_{a=1}^{n-1} \big ( 2 \, (\overline{D}_{e_a} h)(e_a,\overline{\nu}) - (\overline{D}_{\overline{\nu}} h)(e_a,e_a) \big ) \, d\sigma_{\overline{g}} \\ 
&+ \int_{\partial \Omega} \big ( 1 - \frac{1}{2} \, h(\overline{\nu},\overline{\nu}) \big ) \, h(\overline{\nu},\overline{\nu}) \, H_{\overline{g}} \, d\sigma_{\overline{g}} \\ 
&+ \frac{3n-2}{2(n-1)} \int_{\partial \Omega} \sum_{a=1}^{n-1} h(e_a,\overline{\nu})^2 \, H_{\overline{g}} \, d\sigma_{\overline{g}} 
\end{align*} 
by Proposition \ref{boundary.terms.1}. Putting these facts together, the assertion follows. \\

\begin{theorem} 
\label{key.estimate}
Assume that $h$ is divergence-free. Then  
\begin{align*} 
&\bigg | \int_\Omega (R_g - n(n-1)) \, f \, d\text{\rm vol}_{\overline{g}} + \int_{\partial \Omega} (2 - h(\overline{\nu},\overline{\nu})) \, (H_g - H_{\overline{g}}) \, f \, d\sigma \\ 
&+ \frac{1}{4} \int_\Omega |\overline{D} h|^2 \, f \, d\text{\rm vol}_{\overline{g}} + \frac{1}{4} \int_\Omega |\overline{\nabla} (\text{\rm tr}_{\overline{g}}(h))|^2 \, f \, d\text{\rm vol}_{\overline{g}} \\ 
&+ \frac{1}{2} \int_\Omega |h|_{\overline{g}}^2 \, f \, d\text{\rm vol}_{\overline{g}} + \frac{1}{2} \int_\Omega \text{\rm tr}_{\overline{g}}(h)^2 \, f \, d\text{\rm vol}_{\overline{g}} \\ 
&+ \int_{\partial \Omega} h(\overline{\nu},\overline{\nu})^2 \, \partial_{\overline{\nu}} f \, d\sigma_{\overline{g}} + \frac{1}{2} \int_{\partial \Omega} \sum_{a=1}^{n-1} h(e_a,\overline{\nu})^2 \, \partial_\nu f \, d\sigma_{\overline{g}} \\ 
&+ \frac{1}{4} \int_{\partial \Omega} h(\overline{\nu},\overline{\nu})^2 \, H_{\overline{g}} \, f \, d\sigma_{\overline{g}} + \frac{n}{2(n-1)} \int_{\partial \Omega} \sum_{a=1}^{n-1} h(e_a,\overline{\nu})^2 \, H_{\overline{g}} \, f \, d\sigma_{\overline{g}} \bigg | \\ 
&\leq C \int_\Omega |h| \, |\overline{D} h|^2 \, d\text{\rm vol}_{\overline{g}} + C \int_\Omega |h|^3 \, d\text{\rm vol}_{\overline{g}} \\ 
&+ C \int_{\partial \Omega} |h|^2 \, |\overline{D} h| \, d\sigma_{\overline{g}} + C \int_{\partial \Omega} |h|^3 \, d\sigma_{\overline{g}}. 
\end{align*} 
Here, $C$ is a positive constant that depends only on $n$ and $c$.
\end{theorem}

\textbf{Proof.} 
Recall that $f$ is constant along the boundary $\partial \Omega$. Hence, the assertion is a consequence of Proposition \ref{scalar.curvature.3} and Corollary \ref{boundary.terms.2}. \\

\section{Proof of Theorem \ref{main.theorem}}

To prove Theorem \ref{main.theorem}, we need an analogue of Ebin's slice theorem for manifolds with boundary \cite{Ebin} (see also \cite{Fischer-Marsden}). The proof is standard, and works on any compact manifold with boundary. 

\begin{proposition}
\label{slice}
Fix a real number $p > n$. If $\|g - \overline{g}\|_{W^{2,p}(\Omega,\overline{g})}$ is sufficiently small, we can find a diffeomorphism $\varphi: \Omega \to \Omega$ such that $\varphi|_{\partial \Omega} = \text{\rm id}$ and $h = \varphi^*(g) - \overline{g}$ is divergence-free. Moreover, 
\[\|h\|_{W^{2,p}(\Omega,\overline{g})} \leq N \, \|g - \overline{g}\|_{W^{2,p}(\Omega,\overline{g})},\] 
where $N$ is a positive constant that depends only on $\Omega$.
\end{proposition}

\textbf{Proof.} 
Let $\mathscr{S}$ denote the space of symmetric two-tensors on $\Omega$ of class $W^{2,p}$, and let $\mathscr{M}$ denote the space of Riemannian metrics on $\Omega$ of class $W^{2,p}$. Moreover, let $\mathscr{X}$ denote the space of vector fields of class $W^{3,p}$ that vanish along the boundary $\partial \Omega$, and let $\mathscr{D}$ denote the space of all diffeomorphisms $\varphi: \Omega \to \Omega$ of class $W^{3,p}$ satisfying $\varphi|_{\partial \Omega} = \text{\rm id}$. Clearly, the tangent space to $\mathscr{M}$ at $\overline{g}$ can be identified with $\mathscr{S}$; similarly, the tangent space to $\mathscr{D}$ at the identity can be identified with $\mathscr{X}$.

There is a natural action 
\[A: \mathscr{D} \times \mathscr{M} \to \mathscr{M}, \qquad (\varphi,g) \to \varphi^*(g).\] 
Let us consider the linearization of $A$ around the point $(\text{\rm id},\overline{g})$. This gives a map $L: T_{\text{\rm id}} \mathscr{D} \to T_{\overline{g}} \mathscr{M}$. The map $L$ sends a vector field $\xi \in \mathscr{X}$ to the Lie derivative $\mathscr{L}_\xi(\overline{g}) \in \mathscr{S}$. Standard elliptic regularity theory implies that 
\[\mathscr{S} = \{\mathscr{L}_\xi(\overline{g}): \xi \in \mathscr{X}\} \oplus \{h \in \mathscr{S}: \text{\rm $h$ is divergence-free}\}\] 
(compare \cite{Fischer-Marsden}, p.~523). Hence, the assertion follows from the implicit function theorem. \\

We now complete the proof of Theorem \ref{main.theorem}. Let $g$ be a Riemannian metric on the domain $\Omega = \{f \geq c\}$ with the following properties: 
\begin{itemize}
\item $R_g \geq n(n-1)$ at each point in $\Omega$. 
\item $H_g \geq H_{\overline{g}}$ at each point on $\partial \Omega$. 
\item The metrics $g$ and $\overline{g}$ induce the same metric on $\partial \Omega$. 
\end{itemize} 
If $\|g - \overline{g}\|_{W^{2,p}(\Omega,\overline{g})}$ is sufficiently small, Proposition \ref{slice} implies the existence of a diffeomorphism $\varphi: \Omega \to \Omega$ such that $\varphi|_{\partial \Omega} = \text{\rm id}$ and $h = \varphi^*(g) - \overline{g}$ is divergence-free.

Note that $R_{\varphi^*(g)} \geq n(n-1)$ at each point in $\Omega$ and $H_{\varphi^*(g)} \geq H_{\overline{g}}$ at each point on $\partial \Omega$. Applying Theorem \ref{key.estimate} to the metric $\varphi^*(g) = \overline{g} + h$, we obtain 
\begin{align*} 
&\frac{1}{4} \int_\Omega |\overline{D} h|^2 \, f \, d\text{\rm vol}_{\overline{g}} + \frac{1}{4} \int_\Omega |\overline{\nabla} (\text{\rm tr}_{\overline{g}}(h))|^2 \, f \, d\text{\rm vol}_{\overline{g}} \\ 
&+ \frac{1}{2} \int_\Omega |h|_{\overline{g}}^2 \, f \, d\text{\rm vol}_{\overline{g}} + \frac{1}{2} \int_\Omega \text{\rm tr}_{\overline{g}}(h)^2 \, f \, d\text{\rm vol}_{\overline{g}} \\ 
&+ \int_{\partial \Omega} h(\overline{\nu},\overline{\nu})^2 \, \partial_{\overline{\nu}} f \, d\sigma_{\overline{g}} + \frac{1}{2} \int_{\partial \Omega} \sum_{a=1}^{n-1} h(e_a,\overline{\nu})^2 \, \partial_\nu f \, d\sigma_{\overline{g}} \\ 
&+ \frac{1}{4} \int_{\partial \Omega} h(\overline{\nu},\overline{\nu})^2 \, H_{\overline{g}} \, f \, d\sigma_{\overline{g}} + \frac{n}{2(n-1)} \int_{\partial \Omega} \sum_{a=1}^{n-1} h(e_a,\overline{\nu})^2 \, H_{\overline{g}} \, f \, d\sigma_{\overline{g}} \\ 
&\leq C \int_\Omega |h| \, |\overline{D} h|^2 \, d\text{\rm vol}_{\overline{g}} + C \int_\Omega |h|^3 \, d\text{\rm vol}_{\overline{g}} \\ 
&+ C \int_{\partial \Omega} |h|^2 \, |\overline{D} h| \, d\sigma_{\overline{g}} + C \int_{\partial \Omega} |h|^3 \, d\sigma_{\overline{g}}. 
\end{align*} 

If we choose $c \geq \frac{2}{\sqrt{n+3}}$, then 
\[\frac{1}{4} \, H_{\overline{g}} \, f + \partial_{\overline{\nu}} f = \frac{n-1}{4} \, \frac{f^2}{|\overline{\nabla} f|} - |\overline{\nabla} f| = \frac{n-1}{4} \, \frac{c^2}{\sqrt{1-c^2}} - \sqrt{1-c^2} \geq 0\] 
at each point on $\partial \Omega$. This implies 
\begin{align*} 
&\frac{1}{4} \int_\Omega |\overline{D} h|^2 \, f \, d\text{\rm vol}_{\overline{g}} + \frac{1}{4} \int_\Omega |\overline{\nabla} (\text{\rm tr}_{\overline{g}}(h))|^2 \, f \, d\text{\rm vol}_{\overline{g}} \\ 
&+ \frac{1}{2} \int_\Omega |h|_{\overline{g}}^2 \, f \, d\text{\rm vol}_{\overline{g}} + \frac{1}{2} \int_\Omega \text{\rm tr}_{\overline{g}}(h)^2 \, f \, d\text{\rm vol}_{\overline{g}} \\ 
&\leq C \int_\Omega |h| \, |\overline{D} h|^2 \, d\text{\rm vol}_{\overline{g}} + C \int_\Omega |h|^3 \, d\text{\rm vol}_{\overline{g}} \\ 
&+ C \int_{\partial \Omega} |h|^2 \, |\overline{D} h| \, d\sigma_{\overline{g}} + C \int_{\partial \Omega} |h|^3 \, d\sigma_{\overline{g}}. 
\end{align*} 
By the trace theorem, the error terms on the right hand side are bounded from above by $C \, \|h\|_{C^1(\Omega,\overline{g})} \, \|h\|_{W^{1,2}(\Omega,\overline{g})}^2$. Hence, if $\|h\|_{C^1(\Omega,\overline{g})}$ is sufficiently small, then $h$ vanishes identically, and therefore $\varphi^*(g) = \overline{g}$. This completes the proof of Theorem \ref{main.theorem}.

\end{document}